\documentclass[notitlepage,leqno,11pt]{article}
\usepackage{amssymb}
\catcode`\@=11 \@addtoreset{equation}{section}

\catcode`\@=12

\usepackage{latexsym}
\usepackage{textcomp}
\usepackage{amsmath}
\usepackage{amsfonts}
\usepackage{amssymb}
\usepackage{mathrsfs}

\renewcommand{\d}{\delta }

\newcommand{\D }{\Delta }

\newcommand{\e }{\varepsilon }

\renewcommand{\l }{\lambda }

\newcommand{\n }{\nabla }
\newcommand{\vp }{\varphi }

\newcommand{\s }{\sigma }
\newcommand{\Sig }{\Sigma}

\renewcommand{\o }{\omega }
\renewcommand{\O }{\Omega }

\newcommand{\ov}{\overline}

\newcommand{\be}{\begin{equation}}
\newcommand{\ee}{\end{equation}}

\newcommand{\R}{\mathbb{R}}

\newcommand{\de}{\partial}

\newcommand{\ti}{\tilde}

\newcommand{\M}{\mathcal{M}}

\newcommand{\calO }{\mathcal{O}}

\newcommand{\calB }{\mathcal{B}}

\newcommand{\calU}{{\mathcal U}}

\newtheorem{Theorem}{Theorem}[section]
\newtheorem{Lemma}[Theorem]{Lemma}
\newtheorem{Proposition}[Theorem]{Proposition}
\newtheorem{Corollary}[Theorem]{Corollary}
\newtheorem{Remark}[Theorem]{Remark}

\def\proof{\noindent{{\bf Proof. }}}
\def\square{\vbox{
    \hrule height .4pt
    \hbox{\vrule width .4pt height 7pt \kern 7pt
       \vrule width .4pt}
    \hrule height .4pt }}

\def\QED{\hfill {$\square$}\goodbreak \medskip}

\def\R{{\mathbb R}}

\font\sc=cmcsc9  \textwidth=14truecm
\author{Mouhamed Moustapha Fall
\footnote{\footnotesize{Goethe-Universit\"{a}t Frankfurt, Institut f\"{u}r Mathematik.
Robert-Mayer-Str. 10 D-60054 Frankfurt, Germany.
E-mail: \textit{ fall@math.uni-frankfurt.de},
\textit{mouhamed.m.fall@gmail.com}. } }}

\begin{document}

\title{On the Hardy-Poincar\'e inequality with  boundary
singularities }

\date{}
\maketitle

\bigskip

\noindent {\footnotesize{\bf Abstract.}  Let $\O$ be a smooth
bounded domain in $\R^N$ with   $N\ge 1$. In this paper we study the
Hardy-Poincar\'e inequality with  weight function singular at the
boundary of $\O$. In particular we provide sufficient and necessary
conditions on the existence of minimizers.}
\bigskip\bigskip

\section{Introduction}\label{s:i}

Let $\Omega$ be a  domain in $\R^N$, $N\ge 1$, with $0\in\de\O$.
In the framework of Brezis and Marcus \cite{BM}, we  study the
existence and non-existence of  minima for the following
 quotient
\begin{equation}
\label{eq:problem-m} \mu_{\l}(\O):= \inf_{u\in H^{1}_{0}(\O)}
~\frac{\displaystyle\int_{\O}|\nabla u|^2~dx-\l\int_{\O}|u|^2~dx}
{\displaystyle\int_{\O}|x|^{-2}|u|^2~dx}~,
\end{equation}
in terms of   $\lambda\in\R$ and $\O$. The existence and non-existence of extremals
for \eqref{eq:problem-m} were studied in \cite{CaMuPRSE}, \cite{CaMuUMI},
  \cite{Fall}, \cite{FaMu}, \cite{FaMu1}, \cite{NaC}, \cite{Na},   \cite{PT} and the
references there in. Especially in \cite{FaMu}, the authors proved
that for every smooth bounded domain $\O$ of $\R^N$, $N\geq2$, with
$0\in\de\O$
\be\label{eq:supl}
\sup_{\l\in\R}\,\mu_{\l}(\O)=\frac{N^2}{4}=\mu_0\left(\R^N_+\right),
\ee where $\R^N_+=\left\{x\in \R^N\,:\, x^1>0\right\} $, see also
Lemma \ref{lem:Jl1}. In addition they showed that there exists
$\l^*=\l^*(\O)\in[-\infty,+\infty)$ such that
$\mu_\l(\O)<\frac{N^2}{4}$  and it is achieved for all $\l>\l^*$.
 If $\O$ is locally convex at $0$, they proved that
$\l^*\in\R$. Moreover  if $\l^*\in\R$ and $\O$ is
locally concave at $0$ then  there is no minimizer for $\mu_{\l^*}(\O)=\frac{N^2}{4}$.\\
The questions to know whether $\l^*$ is finite for every smooth
domain $\O$ and
 the non-existence of minimizers for $\mu_{\l^*}(\O) $ remained open.\\
We shall show that, indeed, the supremum in \eqref{eq:supl} is
always attained by  $\l^*\in\R$ and that there is no extremals for  $\mu_{\l^*}(\O)$.
Our main result is the following,
\begin{Theorem}\label{th:extJl}
 Let $\O$ be a bounded smooth domain in  $\R^N$, $N\geq2$, with $0\in\de\O$.
 Then there exists $\l^*(\O)\in\R$ such
 that  $\mu_\l(\O)$
 is attained if and only if  $\l>\l^*(\O)$.
\end{Theorem}
We notice that if  $N=1$ then by  \cite{Fall} we have that
$\l^*(\O)\geq0$ and thus $\mu_{\l^*}(\O)$ is not achieved
  by \cite{FaMu1}. We mention that, as observed in \cite{FaMu} and \cite{Fall},
  there are various smooth bounded domains such that $\l^*(\O)<0$. \\
The fact that $\l^*(\O)\in\R$ is a consequence of the following
local Hardy inequality, for $r>0$ small,
\be\label{eq:locHin}
 \int_{\O\cap B_r(0)}|\n
u|^2\,dx\geq\frac{N^2}{4}\int_{\O\cap
B_r(0)}|x|^{-2}|u|^2\,dx\quad\forall \,u\in H^1_0(\O\cap B_r(0)).
\ee
On the other hand the above inequality implies that $\mu_{0}(\O\cap B_r(0))=\frac{N^2}{4}$ by \eqref{eq:supl}.
In particular, even if a domain has negative principal curvatures at
0, its Hardy constant may be equal to $\frac{N^2}{4}$ the  Hardy
constant of the half-space $\R^N_+$. This is not the case for the
Hardy-Sobolev constant, see Ghoussoub-Kang \cite{GK}. Hence the existence of extremals for $\mu_0$ depends on
all the geometry of the domain instead of the  geometric quantities  at the origin, see Proposition \ref{prop:extd}.\\
In Section \ref{s:pn}, we introduce the system of \textit{normal
coordinates} and the modified ground states used in the hall paper.
In section \ref{s:lsReal}, we show that $\l^*\in\R$ and we provide
an improvement of \eqref{eq:locHin}. In Section \ref{s:ne}, we show
that the problem
$$
-\D u-\frac{N^2}{4}|x|^{-2}u=\l u,\quad\textrm{ in } \O
$$
does not possess  a non trivial and nonnegative supersolution in
$H^1_0(\O)\cap C(\O)$. In Section \ref{s:pmt} we prove Theorem
\ref{th:extJl}. Finally in Section \ref{s:Hiww}, we generalize
Theorem \ref{th:extJl} by studying variational problems of type
\eqref{eq:problem-m} with some weights.
%
\section{Preliminaries  and Notations}\label{s:pn}
For $N\geq2$, we denote by  $\{E_1,E_2,\dots,E_N\}$ the standard
orthonormal basis of $\R^N$; $\R^N_+=\{y\in\R^N\,:\,y^1>0\}$;
$B_r(y_0)=\{y\in\R^N\,:\,|y-y_0|<r\}$; $B_r^+=B_r(0)\cap\R^N_+$ and
$S^{N-1}_+=\de B_1(0)\cap \R^N_+$.\\
Let $\calU$ be an open subset of $\R^N$ with boundary $\M:=\de\calU$
a smooth closed hyper-surface of ${\R^N}$ and $0\in\M$. We write
$N_{\M} $ for the unit normal vector-field of $\M$ pointed into
$\calU$. Up to a rotation, we assume that $N_\M(0)=E_1$.  For
$x\in\R^N$, we let $d_\M(x)=\textrm{dist}(\M,x)$ be the distance
function of $\M$. Given $x\in\calU$  and close to $\M$ then it can
be written uniquely as $x=\s_x+d_{\M}(x)\,N_{\M}(\s_x)$, where
$\s_x$ is the projection of $x$ on $\M$. We further   use the
\textrm{Fermi coordinates} $(y^2,\dots,y^N)$ on $\M$ so that for
$\s_x$ close to 0, we have
$$
\s_x=\textrm{Exp}_0\left(\sum_{i=2}^N y^iE_i\right),
$$
where $\textrm{Exp}_0:\R^N\to \M$ is the exponential mapping on $\M$
endowed with the metric induced by $\R^N$. In this way a neighborhood of 0 in $\calU$ can be parameterized by the map
$$
 F_{\M}(y)=\textrm{Exp}_0\left(\sum_{i=2}^Ny^iE_i\right)+y^1\,N_{\M}\left(
\textrm{Exp}_0\left(\sum_{i=2}^Ny^iE_i\right)\right),\quad y\in B_r^+,
$$
for some $r>0$. In this coordinates, the Laplacian $\D$  is given by
$$
\D=\sum_{i=1}^N\frac{\de^2 }{(\de x^i)^2} =\frac{\de^2 }{(\de
y^1)^2}+h_{\M}\circ F_\M\,\frac{\de }{\de
y^1}+\sum_{i,j=2}^N\frac{\de}{\de y^i}\left(\sqrt{|g|}{g^{ij}}
\frac{\de }{\de y^j} \right),
$$
where $h_{\M}(x)=\D\, d_{\M}(x)$; for $i,j=2\dots,N$,
$g_{ij}=\langle\frac{\de F_\M }{\de y^i},\frac{\de F_\M}{\de y^j}
\rangle$; the quantity ${|g|}$ is the determinant of $g$ and
$g^{ij}$ is the component of the
inverse of the matrix $(g_{ij})_{2\leq i,j\leq N}$.\\
Since $g_{ij}=\d_{ij}+O(y^1)+O(|y|^2)$, we have the following
Taylor expansion
\be\label{eq:expDFc}
\D =\sum_{i=1}^N\frac{\de^2 }{(\de y^i)^2}+h_{\M}\circ
F_\M\,\frac{\de }{\de y^1}+\sum_{i=2}^NO_i(|y|)\frac{\de }{\de
y^i}+\sum_{i,j=2}^N O_{ij}(|y|)\frac{\de^2 }{\de y^i\de y^j}.
\ee
%
For $a\in\R$, we put  $X_a(t):=|\log t|^a$, $t\in(0,1)$. Let
$$
\ov{\o}_a(y):=y^1|y|^{-\frac{N}{2}}X_a(|y|)\quad \forall y\in\R^N_+
$$
and put
$$
L_y:=-\sum_{i=1}^N\frac{\de^2 }{(\de y^i)^2}
-\frac{N^2}{4}|y|^{-2}+a(a-1)|y|^{-2}X_{-2}(|y|).
$$
 Then  one easily
verifies that
$$
\begin{cases}
\displaystyle L_y\,\ov{\o}_a=0\,
\textrm{ in }\R^N_+,\\
\displaystyle\ov{\o}_a=0\,\textrm{ on } \de \R^N_+\setminus\{0\},\\
\displaystyle\ov{\o}_a\in H^1(B_R^+)\,\, \forall
R>0,\,\,a<-\frac{1}{2}.
\end{cases}
$$
For $K\in\R$, we define $$\o_{a,K}(y)=e^{K y^1}\,\ov{\o}_a(y).$$
This function satisfies similar boundary and integrability
conditions as $\ov{\o}_a$. In addition  it holds that
%
\be\label{eq:Lyoak}
%
\displaystyle L_y\,\o_{a,K} =-\frac{2K}{y^1}\o_{a,K}
\displaystyle+2K
\left(\frac{N}{2}+aX_{-1}(|y|)\right)\frac{y^1}{|y|^2}\o_{a,K}-K^2\o_{a,K}.
%
\ee
%
Furthermore   for all $a\in\R$
\begin{eqnarray*}
\sum_{i=2}^NO_i(|y|)\frac{\de \o_{a,K}}{\de y^i}+
\sum_{i,j=2}^NO_{ij}(|y|)\frac{\de^2 \o_{a,K}}{\de y^i\de
y^j}&=&y^1\,e^{Ky^1}
O\left(|y|^{-\frac{N}{2}-1}X_a(|y|)\right)\\
&=&\calO_{a,K}(|y|^{-1})\,\o_{a,K}(y).
\end{eqnarray*}
Here the error term $\calO_{a,K}$  has the property that for any
$A>0$, there exist  positive constants $c=c(\O,A,K)$ and
$s_0=s_0(\O,A,K)$ such that
\be\label{eq:errterm}
\left|\calO_{a,K}(s)\right|\leq c\,s\quad\forall
s\in(0,s_0),\,\,\forall a\in[-A,A].
\ee
Let $$W_{a,K}(x):=\o_{a,K}(F_{\M}^{-1}(x)),\qquad\forall\,
x\in\calB_r^+:=F_\M(B_r^+).$$
 Then using \eqref{eq:expDFc}, \eqref{eq:Lyoak} and the
fact that $|x|=|y|+O(|y|^2)$ we obtain the following expansion
%
%
%
%
%
\be\label{eq:expD}
L\,W_{a,K}=-\left(\frac{2K+h_{\M}(x)}{d_{\M}(x)}\right)\,W_{a,K}
+\calO_{a,K}(|x|^{-1})\,W_{a,K}\quad \textrm{ in }\calB_r^+,
\ee
with $ L:=-\D -\frac{N^2}{4}|x|^{-2}+a(a-1)|x|^{-2}X_{-2}(|x|).$
 Moreover it is easy to see that
\be
 \begin{cases}
\label{eq:wH1}
%
W_{a,K}>0\quad\textrm{ in }\calB_r^+,\\
W_{a,K}=0\quad\textrm{ on }\M\cap
\de\calB_r^+\setminus\{0\},\\
W_{a,K}\in H^1(\calB_r^+),\,\,\forall a<-\frac{1}{2}.
\end{cases}
\ee
%
\section{$\l^*(\O)$ is finite}\label{s:lsReal}
 We start with the following local improved Hardy inequality.
\begin{Lemma}\label{lem:loc-hardy}
Let $\calU=\R^N\setminus \ov{B_1(-E_1)}$.
 Then there exist constants $c=c(N)>0$ and
$r_0=r_0(N)>0$ such that for all $r\in(0,r_0)$ the inequality
$$
\int_{\calB_r^+}|\n
u|^2\,dx-\frac{N^2}{4}\int_{\calB_r^+}\frac{|u|^2}{|x|^{2}}\,dx\geq
c\int_{\calB_r^+}\frac{|u|^2}{ |x|^{2}|\log|x||^{2} }\,dx +
(N-1)\,\int_{\calB_r^+}\frac{|u|^2}{d_{\M}(x)}\,dx
$$
holds for all $ u\in H^1_0(\calB_r^+)$.
\end{Lemma}
\proof It is easy to see that $h_{\M}(x)=\frac{N-1}{1+d_\M(x)}$ and
thus
\be\label{eq:ekphgnm1}
-\frac{2(1-N)+h _{\M}(x)}{d_\M(x)}\geq
\frac{N-1}{d_\M(x)}\quad\forall x\in \calU.
\ee
For  $r>0$ small, we set
$$
\ti{w}(x)=\o_{\frac{1}{2},1-N}(F^{-1}_{\M}(x)),\quad \forall
x\in\calB_r^+.
$$
By \eqref{eq:expD} and \eqref{eq:ekphgnm1},  we have
$$
- \frac{\D \ti{w}}{\ti{w}}\geq
\frac{N^2}{4}|x|^{-2}+\frac{1}{4}|x|^{-2}X_{-2}(|x|)+\frac{N-1}{{d_{\M}(x)}}
+O({|x|^{-1}})\,\textrm{ in }\calB_r^+.
$$
 Hence
there exists $r_0=r_0(N)>0$ such that for all $r\in(0,r_0)$
\be\label{eq:dwow}
 - \frac{\D \ti{w}}{\ti{w}}\geq
\frac{N^2}{4}|x|^{-2}+c|x|^{-2}X_{-2}(|x|)+\frac{N-1}{{d_{\M}(x)}}\quad
\textrm{ in }\calB_{r}^+,
\ee
 for some positive constant $c$ depending only on $N$. Fix
$r\in(0,r_0)$ and let $u\in C^\infty_c(\calB_r^+)$. We  put
$\psi=\frac{u}{\ti{w}}$. Then one has $|\n
u|^2=|\ti{w}\n\psi|^2+|\psi\n \ti{w}|^2+\n(\psi^2)\cdot \ti{w} \n
\ti{w}$. Therefore $|\n u|^2=|\ti{w}\n\psi|^2+\n
\ti{w}\cdot\n(\ti{w}\psi^2)$. Integrating by parts, we get
$$
\int_{\calB_r^+}|\n
u|^2\,dx=\int_{\calB_r^+}|\ti{w}\n\psi|^2\,dx+\int_{\calB_r^+}\left(-
\frac{\D \ti{w}}{\ti{w}}\right)u^2\,dx.
$$
 The proof is
then complete by \eqref{eq:dwow} and a desnsity argument.
 \QED
 As a  consequence, we have
\begin{Corollary}\label{cor:loc-hardy-O}
Let $\O$ be Lipschitz domain and of class $C^2$ at $0\in\de\O$.
 Then there exist constants $c=c(\O)>0$ and
$r_0=r_0(\O)>0$ such that for all $r\in(0,r_0)$, the inequality
$$
\int_{\O\cap B_r(0) }|\n
u|^2\,dx-\frac{N^2}{4}\int_{\O\cap B_r(0)}\frac{|u|^2}{|x|^{2}}\,dx\geq
c\int_{\O\cap B_r(0)}\frac{|u|^2}{ |x|^{2}|\log|x||^{2} }\,dx
$$
holds for all $u\in H^1_0(\O\cap B_r(0))$.
\end{Corollary}
\proof Since $\O$ is of class $C^2$ at $0\in\de\O$,  there exits
a ball  with $0\in\de B$ and $\O\subset \calU=\R^N\setminus\ov{B}$.
Therefore by Lemma \ref{lem:loc-hardy}, we get the result. \QED
\begin{Remark}\label{rem:curvature}
We should notice that Lemma \ref{lem:loc-hardy} implies that  "$\O$
is locally concave at $0\in\de\O$" does not necessarly implies that
$\mu(\O)<\frac{N^2}{4}$ as it happens in the Hardy-Sobolev case, see
\cite{GK}, \cite{GR}, \cite{CL}.
\end{Remark}
For sake of completeness, we include the proof of \eqref{eq:supl} in
the  following lemma.
\begin{Lemma}\label{lem:Jl1} Let   $\O$ be a Lipschitz domain and of class $C^2$ at   $0\in\partial\O$.
Then there exists $\l^*(\O)\in\R$ such that
$$
\begin{array}{cc}
\displaystyle  \mu_{\l}(\O)=\frac{N^2}{4}, &  \quad\forall
\l\leq\l^*(\O), \vspace{2mm}\\
 \displaystyle \mu_{\l}(\O)<\frac{N^2}{4}, &
\quad\forall \l>\l^*(\O).
\end{array}
$$
\end{Lemma}

\proof \textbf{Claim:} $\sup_{\l\in\R}\mu_\l\leq \frac{N^2}{4}$.\\
It is well known that $\mu_0(\R^N_+)=\frac{N^2}{4}$, see for
instance \cite{FTT} or \cite{PT}. So for any $\d>0$, we let $u_\d\in
C^\infty_c(\R^N_+)$ such that
$$
\int_{\R^N_+}|\n
u_\d|^2\,dy\leq\left(\frac{N^2}{4}+\d\right)\int_{\R^N_+}|y|^{-2}u_\d^2\,dy.
$$
We let $B$ a ball contained in $\O$ and such that $0\in\de B$.   If
$\e>0$, put
$$
v(x)=\e^{\frac{2-N}{2}}u_\d\left(\e^{-1}F_{\de B}^{-1}(x)\right).
$$
Clearly, provided $\e$ is small enough, we have that $v\in
C^\infty_c(\O)$ thus by the change of variable formula
$$
\mu_\l(\O)\leq\frac{\displaystyle \int_{\O}|\n v|^2\,dx
+\l\int_{\O}v^2\,dx}{\displaystyle \int_{\O}|x|^{-2}v^2\,dx}\leq
\left(1+c\e\right)\frac{\displaystyle \int_{\R^N_+}|\n u_\d|^2\,dy
}{\displaystyle \int_{\R^N_+}|y|^{-2}u_\d^2\,dy }+c\e^2|\l|,
$$
where we have used the fact that $F_{\de B}^{-1}(x)=x+O(|x|^2)$ and
$c$ is a constant depending only on $\O$. We conclude that
$$
\mu_\l(\O)\leq \left(1+c\e\right)
\left(\frac{N^2}{4}+\d\right)+c\e^2|\l|.
$$
Taking the limit in $\e$ and then in $\d$, the claim follows.\\
 \textbf{Claim :} There exists $\ti{\l}\in\R$ such that
 $\mu_{\ti{\l}}=\frac{N^2}{4}$\\
 For $\d>0$ small, we let  $\psi\in C^\infty(B_{\d}(0))$ be a
cut-off function, satisfying
$$
0\leq\psi\leq 1~,\quad\psi\equiv 0~\textrm{in}~\R^N\setminus
B_{\frac{\d}{2}}(0)~,\quad \psi\equiv 1~\textrm{in}~
B_{\frac{\d}{4}}(0)~\!.
$$
We write any $u\in H^1_0(\O)$ as $u=\psi u+(1-\psi)u$, to get
\be\label{eq:uepupumpu} \int_{\O}|x|^{-2}|u|^2~dx\leq
\int_{\O}|x|^{-2}|\psi u|^2~dx+c \int_{\O}|u|^2~dx~\!, \ee
 where the
constant $c$  depends only on $\d$. Since $\psi u \in H^{1}_0(\O\cap
B_\d(0))$,  if  $\d$ is sufficiently small, Corollary
\ref{cor:loc-hardy-O} implies that
\be\label{eq:muCduepupumpu}
  \frac{N^2}{4} \int_{\O}|x|^{-2}|\psi u|^2~dx\leq
\int_{\O}|\nabla(\psi u)|^2~dx.
\ee
 In addition, we have
$$
\int_{\O}|\nabla(\psi u)|^2~dx\leq \int_{\O}|\nabla u|^2~dx+
\frac{1}{2}\int_{\O}\nabla(\psi^2)\cdot\nabla(u^2)~dx+c \int_{\O}
|u|^2~dx~\!.
$$
Using integration by parts we get
$$
\int_{\O}|\nabla(\psi u)|^2~dx\leq \int_{\O}|\nabla
u|^2~dx-\frac{1}{2}\int_{\O}\Delta(\psi^2)|u|^2~dx + c\int_{\O}
|u|^2~dx.
$$
Combining this with (\ref{eq:uepupumpu}) and
(\ref{eq:muCduepupumpu}) we infer that there exits a positive
constant $c$ depending only on $\delta$ and $\O$ such that
$$
  \frac{N^2}{4}\,
{\displaystyle\int_{\O}|x|^{-2}| u|^2~dx}\le {\displaystyle
\int_{\O}|\nabla u|^2~dx+c \displaystyle
\int_{\O}|u|^2~dx}~\!\quad\forall u\in H^1_0(\Omega).
$$
This together with the first calim implies that $\mu_{-c}(\O)=\frac{N^2}{4}$.\\
Finally, noticing that  $\mu_\l(\O)$ is
  decreasing in $\l$, we can set
\be\label{eq:lsdef} \l^*(\O):=\sup\left\{{\l\in\R}\,:\, \mu_\l(\O)=
{\frac{N^2}{4}}\right\}
 \ee
 so that $\mu_\l(\O)<\frac{N^2}{4}$ for all $\l>\l^*(\O)$.
\QED

\section{Non-existence result}\label{s:ne}
In this section we prove the following non-existence result.
\begin{Theorem}\label{th:ne}
Let $\O$ be a bounded Lipschitz  domain of class $C^2$ at
$0\in\de\O$ and let $\l\geq0$.  Suppose that  $u\in H^1_0(\O)\cap
C(\O)$ is a non-negative function
 satisfying
\be\label{eq:ustf}
-\D u-\frac{N^2}{4}|x|^{-2}u\geq-\l u \quad\textrm{in }\O.
\ee
 Then $u\equiv0$.
\end{Theorem}
\proof Up to scaling and rotation, we may assume that $\O$ contains
the ball $B=B_1(E_1)$ such that $\ov{B}\cap\ov{\O}=\{0\}$. We will
use the
coordinates in Section \ref{s:pn} with $\calU=B$ and $\M=\de B$. For  $r>0$ small we define $G_r^+:=F_{\de B}(B_r^+)$.\\
 We suppose  that $u $ does not identically vanish near $0$  and  satisfies
 \eqref{eq:ustf} so that $u>0$ in  $\O\cap B_{r_{0}}(0)$ by  the maximum principle,  for some $r_0>0$.\\
We define
$$
w_{a}(x):=\o_{a,N-1}(F_{\de B}^{-1}(x)),\quad \forall x\in G_r^+.
$$
 Letting
$L:=-\D-\frac{N^2}{4}|x|^{-2}+\l$ then by
 \eqref{eq:expD}
$$
L\,w_a\leq -  \frac{ 2(N-1)+ h_{\de B} }{d_{{\de
B}}}\,{w_a}+\left(\l-\frac{3}{4} \,|x|^{-2}\,X_{-2}(|x|)\right)
\,w_a+\calO_a(|x|^{-1})\,w_a,
$$
for every $a<-\frac{1}{2}$. Since $-h_{\de
B}(x)=(N-1)\left(1+O(|x|)\right)$ in $G_r^+$, by
\eqref{eq:errterm} we can choose $r>0$ small, independent on
$a\in(-1,-\frac{1}{2})$, so that
\be\label{eq:lwaneg}
L\,w_a\leq 0\quad\textrm{ in }G_r^+,\quad\forall
a\in(-1,-\frac{1}{2}).
\ee
 Let $R>0$ so that
$$
R\,{w_a}\leq u\quad\textrm{ on } \ov{F_{\de
B}\left(rS^{N-1}_+\right)}\quad\forall a<-\frac{1}{2}.
$$
By \eqref{eq:wH1}, setting $v_a=R\, {w_a}-u$, it turns out that
$v^+_a=\max(v_a,0)\in H^1_0(G_r^+)$ because $w_a=0$ on $\de B\cap\de
G_r^+$. Moreover by \eqref{eq:ustf} and  \eqref{eq:lwaneg},
$$
L\,v_a\leq 0\quad\textrm{ in }G_r^+,\quad\forall
a\in(-1,-\frac{1}{2}).
$$
 Multiplying the above inequality by $v^+_a$ and integrating by parts
 yields
$$
\int_{G_r^+}|\n
v^+_a|^2\,dx-\frac{N^2}{4}\int_{G_r^+}|x|^{-2}|v^+_a|^2\,dx+\l\int_{G_r^+}|v^+_a|^2\,dx
\leq0.
$$
But then Corollary \ref{cor:loc-hardy-O} implies that $v^+_a=0$ in
$G_r^+$. Therefore $u\geq R\, {w_a}$ for all $a\in(-1,-\frac{1}{2})$
and this contradicts the fact that $\frac{u}{|x|}\in L^2(\O)$
because $\int_{G^+_r}\frac{{w_a}^2}{|x|^{2}}\geq c
\int_{B^+_r}\frac{\o_{a,N-1}^2}{|y|^{2}}\geq\frac{c}{2a+1}|\log
r|^{2a+1}$, for some positive constant $c$ depending only on $B$.
Consequently $u$ vanish identically in $ G_r^+$ and thus by the
maximum principle $u\equiv0$ in  $\O$. \QED
As in \cite{Fall}, starting from exterior domains,
 we can see that, in general, existence of extremals for $\mu_0$ depends  on all the  geometry  of the domain
rather than the geometric constants at the origin. Indeed,
let $G$ be a smooth bounded domain of $\R^N$, $N\geq2$ with $0\in\de G$.
For $r>0$, set $\O_r=B_r(0)\cap  (\R^N\setminus \ov{G})$. It was shown in \cite{Fall} that
there exits $r_1>0$ such that $\mu_0(\O_r)<\frac{N^2}{4}$ for all $r\in(r_1,\infty)$ and $\mu_0(\O_r)$ is achieved.
But  Corollary \ref{cor:loc-hardy-O} and \eqref{eq:supl} yields $\mu_0(\O_r)=\frac{N^2}{4}$ for $r\in(0,r_0)$.
In particular by  Theorem \ref{th:ne}, we get,
\begin{Proposition}\label{prop:extd}
 There exit $r_0,\,r_1>0$ such that the problem
$$
\begin{cases}
\D u+\mu_0(\O_r)\, |x|^{-2}u=0,\quad\textrm{ in } \O_r,\\
u\in H^1_0(\O_r),\\
u\gneqq0\quad\textrm{ in } \O_r
\end{cases}
$$
has a  solution for all $r\in(r_1,\infty)$ and does not have a solution  for every $r\in(0,r_0)$.
\end{Proposition}
\begin{Remark}\label{rem:neGen}
Let $\O$ be as in Theorem \ref{th:ne}. Then by similar argument, one
can show that there is no positive function $u\in H^1_0(\O)\cap
C(\O)$ that satisfies
$$
-\D
u-\frac{N^2}{4}\frac{u}{|x|^{2}}\geq-\frac{\eta(x)}{|x|^{2}}{u}\quad\textrm{
in }\O,
$$
with $\eta$ is continuous, non-negative and
$|\log|x||^{2}\eta(x)\to0$ as $|x|\to0$.
\end{Remark}
\begin{Remark}
We  should mention that some sharp non-existence results of
distributional solution was obtained in \cite{FaMu1}. Indeed assume
that  $\O$ contains a half-ball centered at $0\in\de\O$ and that
$u\in L^2(\O;|x|^{-2}\,dx)$ satisfies
$$
-\int_{\O}u\left(\D \vp+\frac{N^2}{4}\frac{\vp}{|x|^{2}}
\right)\,dx\geq-\frac{3}{4}\int_{\O}u\frac{\vp}{|x|^2|\log|x||^2}\,dx\quad\forall
\vp\in C^\infty_c(\O)
$$
then $u$ vanish in a neighborhood of 0.
\end{Remark}
\section{Proof of  Theorem \ref{th:extJl}}\label{s:pmt}
The proof of the "if" part is similar to the one given in \cite{BM},
see also \cite{FaMu}. Secondly, since the mapping
$\l\mapsto\mu_\l(\O)$ is constant on $(0,\l^*(\O)]$, it is not
difficult to see that $\mu_\l(\O)$ is not achieved for all
$\l<\l^*(\O)$. Now we assume that $\mu_{\l^*}(\O)$ is attained by a
mapping $u\in H^1_0(\O)$. Then it is also achieved by $|u|$ so  we
can assume that $u\gneqq0$. Furthermore since $u$ solves
$$
-\D u - \frac{N^2}{4}|x|^{-2}u=\l^*u\quad\textrm{ in }\O,
$$
by standard elliptic regularity theory, $u$ is smooth in $\O$.
Therefore,  Theorem \ref{th:ne} implies that $u=0$ in $\O$ which is
not possible.\QED
\section{Hardy inequality with weight}\label{s:Hiww}
Let $\O$ be a smooth bounded domain of $\R^N$, $N\geq2$ with
$0\in\de\O$. Following \cite{BM} and \cite{BMS}, we study the
existence of extremals of the following quotient:
\begin{equation}
\label{eq:mpqe} J_{\l}:= \inf_{u\in H^{1}_{0}(\O)}
~\frac{\displaystyle\int_{\O}|\nabla
u|^2p~dx-\l\int_{\O}|x|^{-2}|u|^2\eta~dx}
{\displaystyle\int_{\O}|x|^{-2}|u|^2q~dx}~,
\end{equation}
where the weights $p,q$ and $\eta$ are nonnegative, nontrivial and
satisfy
\begin{equation}\label{eq:ass-pqe}
p\in C^1(\ov{\O}),\quad q,\eta\in C(\ov{\O}), \quad p,\eta>0\textrm{
in }\O\quad\textrm{ and }\quad \eta(0)=0.
\end{equation}
We have the following generalization of Theorem \ref{th:extJl}:
\begin{Theorem}\label{th:mulpqe} Let $\O$ be a smooth bounded
domain of $\R^N$, $N\geq2$ with $0\in\de\O$. Assume that the weight
functions in \eqref{eq:mpqe} satisfy \eqref{eq:ass-pqe}  and that
\begin{equation}\label{eq:pzeqz}
p(0)=q(0)>0.
\end{equation}
Then, there exists $\l^*=\l^*(p,q,\eta,\O)$ such that
$$
\begin{array}{ll}
\displaystyle J_\l=\frac{N^2}{4},\quad\forall\l\leq \l^*,\\
\displaystyle J_\l<\frac{N^2}{4},\quad\forall\l> \l^*.
\end{array}
$$
Furthermore  $J_\l$ is achieved if and only if $\l>\l^*$.
\end{Theorem}
 \proof
\textbf{Step I}: We first show that
\begin{equation}\label{eq:supmlpgelN}
\sup_{\l\in\R}J_\l\leq\frac{N^2}{4}.
\end{equation}
Recall the notation in Section \ref{s:pn}.  For $\rho>0$ small, we
will put $\calB_\rho^+=F_{\de\O}(B_\rho^+)$. By \eqref{eq:pzeqz},
for any $\e>0$ we can let $r_\e>0$ such that
$$
p\leq(1+\e) p(0),\quad q\geq (1-\e)p(0),\quad
\eta\leq\e\quad\textrm{ in } \ov{\calB_{r_\e}^+}.
$$
 By Corollary \ref{cor:loc-hardy-O} and Lemma \ref{lem:Jl1},
$ \mu_0(\calB_{r_\e}^+)= \frac{N^2}{4}$, so for any $\d>0$ we can
let $u\in C^\infty_c(\calB_{r_\e}^+)$ such that
$$
{\int_{\calB_{r_\e}^+}|\n u|^2}{}\leq\left(
\frac{N^2}{4}+\d\right)\int_{\calB_{r_\e}^+}|x|^{-2}u^2.
$$
It turns out that
$$
J_\l\leq \frac{\displaystyle\int_{\O}|\n
u|^2p-\l\int_{\O}|x|^{-2}u^2\eta}{\displaystyle\int_{\O}|x|^{-2}u^2q}\leq
\frac{1+\e}{1-\e}\left(\frac{N^2}{4}+\d
\right)+\frac{\e|\l|}{(1-\e)q(0)}.
$$
Sending $\d$ and  $\e$ to zero, \eqref{eq:supmlpgelN} follows
immediately.\\
 \textbf{Step II}: There exists $\ti{\l}\in\R$ such that
 $J_{\ti{\l}}=\frac{N^2}{4}.$\\
 We fix  $r_0>0$ positive small and put
 \begin{equation}\label{eq:Kzfi}
K_0=\frac{1}{2}\min_{\ov{\calB_{r_0}^+}}\left(- \n p\cdot\n
d_{\de\O}-h_{\de\O}\right).
 \end{equation}
For every $r\in(0,r_0)$, we set
$$
\ti{w}(x)=\o_{\frac{1}{2},K_0}(F_{\de\O}^{-1}(x)),\quad\forall
x\in\calB_r^+.
$$
Notice that $\textrm{div}(p\n \ti{w})=p\D \ti{w}+\n p\cdot\n\ti{w}$.
For $r>0$ small, using \eqref{eq:expD} we get, in $\calB_r^+$,
\begin{equation}\label{eq:div-exp}
-\textrm{div}(p\n\ti{w})=p\frac{N^2}{4}|x|^{-2}\ti{w}+\frac{p}{4}|x|^{-2}X_{-2}(|x|)\ti{w}+\frac{-\n
p\cdot\n
d_{\de\O}-h_{\de\O}-2K_0}{d_{\de\O}}\ti{w}+O(|x|^{-1})\ti{w}.
\end{equation}
Hence by \eqref{eq:pzeqz} and \eqref{eq:Kzfi}  there exist constants
$c>0$ and $r_1>0$ (depending on $p$, $q$, $\eta$ and $\O$) such that
for all $r\in(0,r_1)$
\begin{equation}\label{eq:divptwgeq}
-\textrm{div}(p\n\ti{w})\geq
q\frac{N^2}{4}|x|^{-2}\ti{w}+c|x|^{-2}X_{-2}(|x|)\ti{w}\quad
\calB_r^+.
\end{equation}
Fix $r\in(0,r_1)$ and let $u\in C^\infty_c(\calB_r^+)$. We  put
$\psi=\frac{u}{\ti{w}}$. Then one has $|\n
u|^2=|\ti{w}\n\psi|^2+|\psi\n \ti{w}|^2+\n(\psi^2)\cdot \ti{w} \n
\ti{w}$. Therefore $|\n u|^2p=|\ti{w}\n\psi|^2p+p\n
\ti{w}\cdot\n(\ti{w}\psi^2)$. Integrating by parts, we get
$$
\int_{\calB_r^+}|\n
u|^2p\,dx=\int_{\calB_r^+}|\ti{w}\n\psi|^2p\,dx+\int_{\calB_r^+}\left(-
\frac{\textrm{div}(p\n \ti{w})}{\ti{w}}\right)u^2\,dx.
$$
This together with \eqref{eq:divptwgeq} yields
\begin{equation}\label{eq:locHG}
\int_{\calB_r^+}|\n u|^2p\,dx\geq
\frac{N^2}{4}\int_{\calB_r^+}|x|^{-2}u^2q\,dx+c
\int_{\calB_r^+}|x|^{-2}X_{-2}(|x|)u^2.
\end{equation}
We can now proceed as in the proof of Lemma \ref{lem:Jl1} (since
$\eta>0$ in ${\O}$) to conclude that there exists a constant
$C=C(p,q,\eta,\O)>0$ such that
$$
  \frac{N^2}{4}\,
{\displaystyle\int_{\O}|x|^{-2} u^2q~dx}\le {\displaystyle
\int_{\O}|\nabla u|^2p~dx+C \displaystyle
\int_{\O}|x|^{-2}u^2\eta~dx}~\!\quad\forall u\in H^1_0(\Omega).
$$
Therefore we can define $\l^*$ as in \eqref{eq:lsdef} to end the
proof of this step.\\
\textbf{Step III:} Let $u\in H^1_0(\O)\cap C(\O)$ is a non-negative
function
 satisfying
\be\label{eq:ustfpqe}
-\textrm{div}( p\n u)-\frac{N^2}{4}q|x|^{-2}u\geq-\l|x|^{-2}\eta u
\quad\textrm{in }\O.
\ee
 Then $u\equiv0$.\\
Here, we assume that $\O$ contains the ball $B=B_1(E_1)$ such that
$\ov{B}\cap\ov{\O}=\{0\}$ and set $G_r^+=F_{\de B}(B_r^+)$. As in
the previous step, we put
 \begin{equation}\label{eq:Kofi}
K_1=\frac{1}{2}\max_{\ov{G_{r_0}^+}}\left(- \n p\cdot\n
d_{\de\O}-h_{\de\O}\right).
 \end{equation}
 For $r\in(0,r_0)$  and
$a<-\frac{1}{2}$, we set
$$
{w_a}(x)=\o_{a,K_1}(F_{\de B}^{-1}(x)),\quad\forall x\in
G_r^+=F_{\de B}(B_r^+).
$$
Letting $L=-\textrm{div}(p\n\cdot
)-\frac{N^2}{4}q|x|^{-2}+|\l||x|^{-2}\eta$ then by \eqref{eq:Kofi}
and  \eqref{eq:pzeqz}, we get
$$
Lw_a\leq\left(|\l||x|^{-2}\eta-\frac{3}{4}p|x|^{-2}X_{-2}(|x|)
\right)w_a+\calO_a(|x|^{-1})w_a\quad\textrm{ in } G_r^+.
$$
Therefore  by \eqref{eq:errterm} we can choose $r>0$ small,
independent on $a\in(-1,-\frac{1}{2})$, so that
\be\label{eq:lwanegpqe}
L\,w_a\leq 0\quad\textrm{ in }G_r^+,\quad\forall
a\in(-1,-\frac{1}{2}).
\ee
If  $u\gvertneqq0$ near the origin then by the maximum principle, we
can assume that $u>0$ in $G_{2r}^+$. Hence we can let $R>0$ so that
$$
R\,{w_a}\leq u\quad\textrm{ on } \ov{F_{\de
B}\left(rS^{N-1}_+\right)}\quad\forall a<-\frac{1}{2}.
$$
By \eqref{eq:wH1}, setting $v_a=R\, {w_a}-u$, it turns out that
$v^+_a=\max(v_a,0)\in H^1_0(G_r^+)$. Moreover by \eqref{eq:ustfpqe}
and  \eqref{eq:lwanegpqe},
$$
L\,v_a\leq 0\quad\textrm{ in }G_r^+,\quad\forall
a\in(-1,-\frac{1}{2}).
$$
 Multiplying the above inequality by $v^+_a$ and integrating by parts
 yields
$$
\int_{G_r^+}|\n
v^+_a|^2p\,dx-\frac{N^2}{4}\int_{G_r^+}|x|^{-2}|v^+_a|^2q\,dx+|\l|\int_{G_r^+}|x|^{-2}|v^+_a|^2\eta\,dx
\leq0.
$$
But then \eqref{eq:locHG} implies that $v^+_a=0$ in $G_r^+$.
Therefore $u\geq R\, {w_a}$ for all $a\in(-1,-\frac{1}{2})$ and this
contradicts the fact that $\frac{u}{|x|}\in L^2(\O)$. Consequently
$u$ vanish identically in $ G_r^+$ and thus by the
maximum principle $u\equiv0$ in  $\O$.\\
\textbf{Step IV:} If $J_\l<\frac{N^2}{4}$ then it is achieved.\\
The proof of the existence part, since $\eta(0)=0$, is similar to
the one given in \cite{BM} so we skip it.
 \QED
\begin{Remark}
Let $\O$ be a smooth smooth bounded domain of $\R^N$, $N\geq2$. Let
$\Sig_k$ be a smooth compact sub-manifold of $\de\O$ with dimension
$0\leq k\leq N-1$. Here $\Sig_0$ is a single point. Consider the
problem $(P_k^\l)$ of finding minimizers for the quotient:
\begin{equation}
\label{eq:mpqek} J_{\l}^k:= \inf_{u\in H^{1}_{0}(\O)}
~\frac{\displaystyle\int_{\O}|\nabla
u|^2p~dx-\l\int_{\O}\textrm{dist}(x,\Sig_k)^{-2}|u|^2\eta~dx}
{\displaystyle\int_{\O}\textrm{dist}(x,\Sig_k)^{-2}|u|^2q~dx}~,
\end{equation}
where the weights $p,q$ and $\eta$ are smooth positive in $\ov{\O}$
with $\eta=0$ on $\Sig_k$ and the following  normalization
\begin{equation}
\min_{\Sig_k}\frac{p}{q}=1
\end{equation}
holds. We put
\begin{equation}
I_{k}=\int_{\Sig_k}\frac{d\s}{\sqrt{1-\left(q(\s)/p(\s)\right)}},\quad
1\leq k\leq N-1\quad\textrm{ and }\quad I_0=\infty.
\end{equation}
It was shown  in \cite{BM} that there exists $\l^*$ such that
 if $\l>\l^*$ then
$J^{N-1}_\l<\frac{1}{4}$ and $(P_{N-1}^\l)$ has a solution while for
$\l\leq\l^*$, $J^{N-1}_\l=\frac{1}{4}$ and  $(P_{N-1}^\l)$ does not
have a solution whenever $\l<\l^*$. The critical case
$(P_{N-1}^{\l^*})$ was treated in \cite{BMS}, where the authors
proved that $(P_{N-1}^{\l^*})$ admits a solution if and only if
$I_{N-1}<\infty$. This clearly holds here for $(P_{0}^{\l^*})$  by
Theorem \ref{th:mulpqe}. We believe that such type of results remain
true for all $k$  by taking in to account that in the flat case,
$$
\inf_{u\in H^1_0(\R^N_+)}\frac{\displaystyle\int_{\R^N_+}|\n
u|^2\,dx}{\displaystyle\int_{\R^N_+}\frac{u^2}{x_1^2+\cdots+x_{N-k}^2}\,dx}=\frac{(N-k)^2}{4},
$$
see \cite{FTT}, with  $\R^N_+=\left\{x\in \R^N\,:\, x^1>0\right\} $.
\end{Remark}

\textbf{ Acknowledgments} \\
The author would like to thank Professor Haim Brezis
 for his comments and suggestions. 
 He is grateful to the referee for his comments.  This work is partially supported by the Alexander-von-Humboldt Foundation.


\begin{thebibliography}
\footnotesize
%
%
\bibitem{BM} Brezis H. and Marcus M., Hardy's inequalities revisited.
 Dedicated to Ennio De Giorgi.
  Ann. Scuola Norm. Sup. Pisa Cl. Sci. (4)  25  (1997),  no. 1-2, 217-237.
  \bibitem{BMS}  Brezis H.,  Marcus M. and  Shafrir I.,
   Extermal functions for Hardy's inequality with weight, J. Funct. Anal. 171 (2000), 177-191.
%

%
%
%
\bibitem{CaMuPRSE}
{  Caldiroli P.,  Musina R.}, On a class of 2-dimensional singular
elliptic problems. Proc. Roy. Soc. Edinburgh Sect. A 131 (2001),
479-497.
%

\bibitem{CaMuUMI} Caldiroli P., Musina R., Stationary states for a
two-dimensional~ singular Schr\"odinger equation. Boll. Unione Mat.
Ital. Sez. B Artic. Ric. Mat.  (8) 4-B (2001), 609-633.
%
\bibitem{CL} Chern J-L. and   Lin C-S., Minimizers of Caffarelli-Kohn-Nirenberg Inequalities with the
Singularity on the Boundary. Archive for Rational Mechanics and
Analysis Volume 197, Number 2 (2010), 401-432.
%
%
%
\bibitem{Fall} Fall M. M., A note on Hardy's inequalities with
boundary singularities. Pr\'epublication D\'epartement de
Math\'ematique Universit\'e Catholique de Louvain-La-Neuve 365
(2010), {\tt http://www.uclouvain.be/38324.html}.
%
\bibitem{FaMu} Fall M. M., Musina R., Hardy-Poincar\'e inequalities with boundary singularities.
Pr\'epublication D\'epartement de Math\'ematique Universit\'e
Catholique de Louvain-La-Neuve 364 (2010), {\tt
http://www.uclouvain.be/38324.html}.
%
\bibitem{FaMu1} Fall M. M., Musina R., Sharp nonexistence results for a linear elliptic inequality involving Hardy and Leray potentials.
Journal of Inequalities and Applications, vol. 2011, Article ID 917201, 21 pages, 2011.
 doi:10.1155/2011/917201.
%
%
\bibitem{FTT}   Filippas S.;  Tertikas A. and
Tidblom J., On the structure of Hardy-Sobolev-Maz'ya inequalities .
J. Eur. Math. Soc., 11(6), (2009), 1165-1185.
%
%
%
%
\bibitem{GK} Ghoussoub N., Kang X. S., Hardy-Sobolev critical elliptic equations with boundary singularities.
Ann. Inst. H. Poincar\'e Anal. Non Lin\'eaire  21  (2004),  no. 6,
767--793.
%
\bibitem{GR} Ghoussoub N. and Robert F.: The effect of curvature on the best
constant in the Hardy– Sobolev inequalities. Geom. Funct. Anal.
16(6), 1201-1245 (2006).
%
%
%
%
%
%
\bibitem{NaC}  Nazarov A. I., Hardy-Sobolev Inequalities in a cone, J. Math.
Sciences, 132, (2006), (4), 419-427.
%
\bibitem{Na} Nazarov A. I., Dirichlet and Neumann problems to critical
Emden-Fowler type equations. J Glob Optim (2008) 40, 289-303.
%
\bibitem{PT}{Pinchover Y.,  Tintarev K.},
 { Existence of minimizers for Schr\"{o}dinger operators under domain perturbations with application
  to Hardy's inequality}.
 Indiana Univ. Math. J. 54 (2005), 1061-1074.
%

%
%
%
%
%
%
%
%
%
%
%
%
\end{thebibliography}
\end{document}